\documentclass[a4paper,12pt]{article}
\usepackage{amsmath,amssymb}
\usepackage{graphics}

\def\R{\mathbb R}
\def\C{\mathbb C}

\def\qed{
\hspace{7mm} {\lower 1.5pt \hbox{\scalebox{.65}[1.4]{$\blacksquare$}}}}

\begin{document}
\begin{center}
{\large Boots-Strapping Argument for a $J$-Holomorphic Mapping on a Strongly Pseudo-Convex Manifold}
\end{center}

\vspace{3mm}
\begin{center}
{\large Takanari Saotome}\\
Institute of Mathematics, University of Tsukuba\\
1-1-1, Tennoudai, Tsukuba, Ibaraki, 305-8571.
\end{center}

\vspace{5mm}\noindent
{\bf Abstract.}\hspace{2mm} For $J$-holomorphic mappings for a strongly pseudo-convex manifold, we prove elliptic regularity by the argument of boots-strapping.

\vspace{2mm}\indent
{\it key words:}\hspace{2mm} pseudo-convex manifold; pseudo-holomorphic

\renewcommand{\thefootnote}{}
\footnote{
\hspace{-6.5mm}2000 Mathematics Subject Classification : 32V05; 53D35.\\
tsaotome@math.tsukuba.ac.jp
}

\vspace{5mm}\noindent
{\bf \S 0.\hspace{2mm}Introduction.}
\setcounter{equation}{0}
\renewcommand{\theequation}{0-\arabic{equation}}

\vspace{2mm}
The theory of $J$-holomorphic curves is one of the most developed subject in the study of symplectic geometry ever since its inception in Gromov's paper of 1985 [6].
This celebrated theory is based on meny analytic properties of $J$-holomorphic curves.
In perticular, the removable singularity theorem and the Fredholm property of the equation of $J$-holomorphic curves are one of the most important facts to see the moduli space is acutually a compact smooth finite dimensional manifold.

Similarly to symplectic manifolds, strongly pseudo-convex manifolds admit a "nondegenerate clsed $2$-form" and an "almost complex structure".
So the definition of $J$-holomorphy makes sence, and we can see that $J$-holomorphic mappings for strongly pseudo-convex manifolds also have useful properties similar to the symplectic case [11].
This paper is one of the result to construct an analogous theory in odd-dimensional geometry, we consider a $J$-holomorphic mapping from a $3$-dimensional Sasakian manifold into a strongly pseudo-convex manifold.
In my previous paper [12], we have seen that the removable singularity theorem hold for $J$-holomorphic mappings for strongly pseudo-convex manifolds.

\vspace{2mm}
It also can be seen that, similarly to $J$-holomorphic curves in symplectic manifold, $J$-holomorphic mappings satisfy an elliptic differential equation and so have hypo-ellipticity.
Although the hypo-ellipticic regularity assures the smoothness of solutions of a partial differential equation, the assuption of this regularity requires "weak smoothness" of the solutions.
In the even-dimensional geometry, the argument of boots-strapping appears to avoid this difficulty.

In this paper, we will see the boots-strapping argument also valid in the odd-dimensional geometry.
More explicitely, we will obtain following results. (For precise definitions, see \S 2.)

\vspace{3mm}\noindent
{\bf T{\sc heorem}\hspace{1mm}1}\hspace{2mm}Fix an integer $l \ge 2$ and a number $p > 3$.
Let $J \in \mathcal{J}^l(M)$ and suppose $u : \Sigma \longrightarrow M$ is a $W^{1, p}$-mapping such that  
\begin{enumerate}
 \item \ $du \circ j = J \circ du,$ on $P_\Sigma$\hfill(jh-1)
 \item \ $du(\xi_\Sigma) = \lambda \xi,$ for a positive smooth function $\lambda > 0$ on $\Sigma$\hfill(jh-2)
 \item \ $u^\ast \theta = \lambda \theta_\Sigma$.\hfill(jh-3)
\end{enumerate}
Then $u$ is of class $W^{l, p}.$

\vspace{5mm}
The author would like to thank professors Mitsuhiro Itoh for pertinent comments and for encouraging steadily the author of this article.

\vspace{5mm}
\noindent
{\bf \S 1.\hspace{2mm}Preliminary.}
\setcounter{equation}{0}
\renewcommand{\theequation}{1-\arabic{equation}}

\vspace{2mm}
We first recall basic notions of a strongly pseudo-convex manifold $M$,

A $(2n + 1)$-dimensional {\it strongly pseudo-convex manifold} $M$ is an oriented smooth manifold which carries a strongly pseudo-convex structure $(P, J, \theta)$.
We call a triple $(P, J, \theta)$ of a contact $1$-form and an almost complex strucutre $J$ on a distribution $P = \ker{\theta}$ strongly pseudo-convex strucutre if $J$ satisfies the integrable condition
\[ [X, Y] - [JX, JY] - J([JX, Y] + [X, JY]) = 0,\ \mbox{for}\ X, Y \in \Gamma(P),\]
and the {\it Levi-form} $L(X, Y) = -d\theta(JX, Y)$ for $X,Y \in P$, is positive definite.
We call, in what follows, a strongly pseudo-convex manifold s.p.c. manifold in abbreviation.

Consider the complexification of $TM$ and its $\sqrt{-1}$ eigen-subspace $S = \{ X - \sqrt{-1}JX \mid X \in P\} \subset \mathbb{C}TM$ .
It holds $S \cap \overline{S} = \{ 0\}$ and $[\Gamma(S), \Gamma(S)] \subset \Gamma(S)$, where $\overline{S}$ is the complex conjugation of $S$.

If $(M, \theta, J)$ is an s.p.c. manifold, then there exists a unique nonvanishing vector field $\xi$ which satisfies $\theta(\xi) = 1$ and $i(\xi)\, d\theta = 0$.
We call this vector field the characteristic vector field.
(The characteristic vector field is sometimes called Reeb vector field.)
By the definition, we have a decomposition $TM = P \oplus \R \xi$, and hereafter we consider the Levi-form as a tensor on $TM$ by extending $J$ on $TM$ as $J\xi = 0.$
Then, we have a canonical Riemannian metric $g = g_{(\theta, J)}$, defined by
\[ g(X, Y) = L(X, Y) + \theta(X) \theta(Y).\]

Since the contact form $\theta$ and the characteristic vector field $\xi$ satisfy $L_\xi \theta = 0,$ we can see that the necessary and sufficient condition of $L_\xi J = 0$ is $\xi$ is Killing with respect to the metric $g.$
We call an almost complex structure $J$ {\it normal} or {\it K-contact} when $J$ and $\theta$ define the metric $g$ for which $\xi$ is Killing.
We call an s.p.c. manifold with normal almost complex structure {\it Sasakian manifold}.

\vspace{5mm}
Let $(\Sigma^3, \theta_\Sigma, j)$ be a connected $3$-dimensional Sasakian manifold and $(M, \theta, J)$ be an s.p.c. manifold.
Consider a smooth mapping $u : \Sigma \longrightarrow M$ which obeys the conditions
\begin{enumerate}
 \item \ $du \circ j = J \circ du,$ on $P_\Sigma$\hfill(jh-1)
 \item \ $du(\xi_\Sigma) = \lambda \xi,$ for a positive smooth function $\lambda > 0$ on $\Sigma$\hfill(jh-2)
 \item \ $u^\ast \theta = \lambda \theta_\Sigma$.\hfill(jh-3)
\end{enumerate}
Here, we note that the equation (jh-1) is just an analogue of the notion of $J$-holomorphic curve in symplectic geometry, which M. Gromov initially defined in [7].
So, we call such maps {\it $J$-holomorphic mapping for an s.p.c. manifold $M$}.

By the conditions (jh-1) and (jh-3), $J$-holomorphic mappings preserve the holomorphic structures, i.e. we have $u_\ast(S_\Sigma) \subset S$ for a $J$-holomorphic mapping $u$.
This condition is often referred as CR-holomorphy of a mapping between CR manifolds.

\vspace{2mm}
\noindent
{\bf Remark.}\hspace{1mm}Here we note that if the function $\lambda$ is identically $1$ over $\Sigma$, then the $J$-holomorphic mapping $u : \Sigma \longrightarrow M$ is an isometric immersion.

\vspace{2mm}
\noindent{\bf Examples.}

i)\hspace{2pt}{\it Plane sections}:\hspace{2mm}Let $f(z_1, \cdots , z_{n + 1})$ be a weighted homogeneous polynomial with an isolated singular point $0 \in \C^{n + 1}$.
It is well known that the link of the locus $M^\prime = f^{-1}(0)$ has an s.p.c. structure.

Let $M = S^{2n + 1} \cap M^\prime$ be the link and $u : \Sigma \longrightarrow M$ be a natural inclusion map, where $\Sigma$ is the link of the locus of the weighted homogeneous polynomial $f_0(z_1, z_2, z_3) = f(z_1, z_2, z_3, 0, \cdots , 0)$, the restriction of $f$ to a four dimensional subspace $V$.
So $\Sigma$ is a section of $M$ by $V$.
Then it is obvious that $u$ is a $J$-holomorphic mapping into $M$.

\vspace{2mm}
ii)\hspace{2pt}{\it Veronese mappings}:\hspace{2mm}Let $\Sigma = S^3 \subset \C^2$ and $M = S^5 \subset \C^3$ be standerd spheres.
Then we can see the Veronese mapping $v : S^3 \longrightarrow S^5$
\[ v(z_1, z_2) = (z_1^2, \sqrt{2}z_1z_2, z_2^2),\]
satisfies the condition of $J$-holomorphic mappings, by straight computation.
Then it easily can be shown this mapping preserves the norm of complex Eucilidian spaces and the holomorphic strucutre.
It can be also seen that $dv(\xi_\Sigma) = 2\xi$ so that $v$ is a $J$-holomorphic mapping between s.p.c. manifolds. 

More generally, the lift of $J$-holomorphic curves can be considered as an example of $J$-holomorphic mappings.
Let $\Sigma^3 \longrightarrow \underline{\Sigma}^2$ and $M^{2n + 1} \longrightarrow \underline{M}^{2n}$ be negative $S^1$ bundles over a Riemann surface and a $2n$-dimensional K\"{a}hler manifold, respectively.
Let $\underline{u} : \underline{\Sigma} \longrightarrow \underline{M}$ be a holomorphic curve between base spaces and we assume that there exists a bundle isomorphism from the pull-back bundle $\underline{u}^\ast M$ to the bundle $\Sigma$ which preserves connections.
Then we can see that the bundle isomorphism induces a lift $u : \Sigma \longrightarrow \underline{u}^\ast M \subset M$ over the holomorphic curve $\underline{u},$ which is $J$-holomorphic. 

\vspace{5mm}
For $J$-holomorphic mappings between s.p.c. manifolds, we have simple but strong property as follows:
Namely, the condition (jh-2) and (jh-3) imply that the function $\lambda$ must be a constant.

\vspace{2mm}
\noindent
{\bf P{\sc roposition}\hspace{1mm}1.1}. For every mapping $u : \Sigma \longrightarrow M$ satisfying (jh-2) and (jh-3), the function $\lambda$ defined above is constant.

\vspace{1mm}
\noindent
{\bf Proof.}\hspace{2mm} In fact, for $X \in P_\Sigma$ we have $X \lambda = 0$, since
\begin{eqnarray*}
0 &=& d\theta(du (\xi_\Sigma), du X) = u^\ast d\theta(\xi_\Sigma, X)\\
&=& (d \lambda \wedge \theta_\Sigma + \lambda d\theta_\Sigma)(\xi_\Sigma, X) = X \lambda,
\end{eqnarray*}
Therefore by the strong pseudo-convexity of $\Sigma$ we deduce $\xi \lambda = 0,$ and hence $d\lambda = 0.$\qed

\vspace{5mm}
\noindent
{\bf \S 2.\hspace{2mm}Proof of Theorem 1.}
\setcounter{equation}{0}
\renewcommand{\theequation}{3-\arabic{equation}}

\vspace{2mm}
In this section we shall prove the regularity property of $J$-holomorphic mappings, that is useful to prove some kind of compactness theorem of the moduli space of $J$-mappings.
We mainly follow the argument which used in the symplectic geometry, and the basic idea of the proof is that to reduce the problem to regularity of the usual laplace equation.

In what follows, we will concentrate to local situation.
More precisely, we consider as
\[ u : \Omega \longrightarrow \R^{2n + 1},\quad J \in W^{1, p}(\Omega; \R^{(2n + 1) \times (2n + 1)}),\]
where $\Omega \subset \R^3$ is a connected open set.
We note that in this situation we can take a coordinate $(s, t, z),\ (x^1, \ldots , x^n, y^1, \ldots , y^n, w)$ such that $\partial/\partial z = \xi_\Sigma,\ \partial/\partial z = \xi.$

We first show the following two technical lemmas.

\vspace{2mm}
\noindent
\textbf{Lemma 2.1}\hspace{2mm}
For given $p > 3,$ there exists a finite sequence $\{ q_i\}_{i = 0,\ \ldots \ , m}$ such that
\[ q_i < q_{i + 1},\quad \frac{p}{p - 1} < q_0 \le p,\quad q_{m - 1} < \frac{3p}{p - 3} < q_m.\]

\vspace{1mm}
\noindent
{\bf Proof.}\hspace{2mm}
We note that for any $p, q, r$, we have
\[ \frac{1}{p} + \frac{1}{q} = \frac{1}{r} \iff r = \frac{pq}{p + q},\]
and for such $p, q, r$ and a positive integer $n < p,$ 
\[ r > n \iff q > \frac{np}{p - n}.\]
Therefore, we can take $q_0$ so that $\displaystyle \frac{p}{p - 1} < q_0 \le p$.
Define a sequence $\{ q_i\}$ by
\[ q_{i + 1} = \frac{3r_i}{3 - r_i},\quad r_i = \frac{pq_I}{p + q_i}.\]

Now since $p > 3$, the function $a : \R \longrightarrow \R$ defined by $\displaystyle a(q) = \frac{3r}{3 - r},\ r = \frac{pq}{p + q}$, satisfies $a(q) > q$ for $q > 0$.
Moreover, we have $\displaystyle q_2 < q_1 < \frac{3p}{p - 3} \ \Rightarrow \ a(q_2) < a(q_1).$
This facts imply there exists a finite number $m$ such that $\displaystyle \frac{3p}{p - 3} < q_m,$ and this completes the proof. \qed

\vspace{2mm}
\noindent
\textbf{Lemma 2.2}\hspace{2mm}
Let $p > 3$ and $1 < r \le p.$
Then for $f \in W^{1, p}(\R^3)$ and $g \in W^{1, r}(\R^3)$, we have $fg \in W^{1, r}(\R^3)$.
Moreover, there exists a constant $c$, depend on $p, r$, such that
\[ \| fg\|_{W^{1, r}} \le c\| f\|_{W^{1, p}} \| g\|_{W^{1, r}}.\]

\vspace{1mm}
\noindent
{\bf Proof.}\hspace{2mm}

\vspace{2mm}
To prove Theorem 1, we will use some classical analytic theorems. 
We state them without proofs for the sake of readers. (see for more details Appendix B of [10].)

\vspace{2mm}
\noindent
\textbf{Theorem 2.3}\hspace{2mm}
Let $\Omega \subset \R^n$ be a bounded Lipschitz domain and suppose that $kp < n.$
Then there exists a constant $c > 0$, depend only on $k, p$ and $\Omega,$ such that
\[ \| u\|_{L^{np/(n - kp)}} \le c\| u\|_{W^{k, p}}\]
for $u \in C^\infty(\overline{\Omega}).$
If $\displaystyle q < \frac{np}{n - kp}$ then the inclusion $W^{k, p}(\Omega) \longrightarrow L^q(\Omega)$ is compact.

These are called Sobolev estimate and Rellich's embedding theorem.

\vspace{2mm}
\noindent
\textbf{Theorem 2.4}\hspace{2mm}
Let $1 < r < \infty$ and $\Omega^\prime \subset \Omega \subset \R^n$ be open sets such that $\overline{\Omega^\prime} \subset \Omega.$
Then there exists a constant $c > 0$ such that the following holds.
Assume that $u \in L^1_{loc}(\Omega)$ and $f = (f_0, \ldots \ , f_n) \in L^p_{loc}(\Omega, \R^{n + 1})$ satisfy
\[ \int_\Omega u(x) \triangle \varphi (x)\, dx = \int_\Omega f_0(x) \varphi (x)\, dx - \sum^n_{j = 1} \int_\Omega f_j(x) \partial_j \varphi (x)\, dx\]
for every $\varphi \in C^\infty_0(\Omega).$
Then $u \in W^{1, p}_{loc}(\Omega)$ and
\[ \| u\|_{W^{1, p}(\Omega^\prime)} \le c( \| f\|_{L^p(\Omega)} + \| u\|_{L^p(\Omega)} ).\]

\vspace{2mm}
Next we show the regularity of "weak $J$-holomorphic mappings".
The proof of theorem 1 is based on this key-lemma.

\vspace{2mm}
\noindent
\textbf{Proposition 2.5}\hspace{2mm}
Let $p > 3$ and $\infty > r > 1,\ \displaystyle \frac{1}{p} + \frac{1}{q} = \frac{1}{r}.$
We assume $J \in W^{1, p}(\Omega; \R^{(2n + 1) \times (2n + 1)}), u \in L^q_{loc}(\Omega, \R^{2n + 1})$ and there exist $\eta, \nu \in L^r_{loc}(\Omega, \R^{2n + 1})$ such that
\begin{eqnarray*}
&& \int_\Omega \langle \partial_s \varphi + J^T \partial_t \varphi, u \rangle = -\int_\Omega \langle \varphi, \eta + (\partial_t J)u \rangle,\\
&& \int_\Omega \langle \partial_t ((J^T)^2 + 1) \varphi, u \rangle = -\int_\Omega \langle \varphi, \nu + \theta_\Sigma(\partial_t) \partial_z u \rangle,\quad \mbox{and}\\
&& \partial_z u = \mbox{const}.
\end{eqnarray*}
for any $\varphi \in C^\infty_0(\Omega, \R^{2n + 1}).$
Then we have $u \in W^{1, r}_{loc}(\Omega, \R^{2n + 1}).$

\vspace{1mm}
\noindent
{\bf Proof.}\hspace{2mm}
According to the theorem , it is suffice to examine $\displaystyle \int_\Omega \langle \triangle \psi, u \rangle$.
Let $\varphi = \partial_s \psi + J^T \partial_t \psi$.
Then 
\[ \partial_s \varphi + J^T \partial_t \varphi = \square \psi - (\partial_s J^T) \partial_t \psi + (\partial_t J^T) J^T \partial_t \psi - \partial_t K \partial_t \psi, \]
where $\square = (\partial_s)^2 + (\partial_t)^2$ is a sub-Laplacian and $K = (J^T)^2 + 1.$

By the hypothesis of $\partial_z u = \mbox{const}$, we note that $\displaystyle \int_\Omega \langle \square \psi, u \rangle = \int_\Omega \langle \triangle \psi, u \rangle = \int_\Omega \langle  [(\partial_s)^2 + (\partial_t)^2 + (\partial_z)^2] \psi, u \rangle.$
So applying the condition of $K$ to the right hand side of previous formula
\begin{eqnarray*}
\int_\Omega \langle \triangle \psi, u \rangle &=& -\int_\Omega \langle \partial_s \psi, (\partial_t J) u + \eta \rangle\\
&&\hspace{20pt}  - \int_\Omega \langle \partial_t \psi, -(\partial_s J) u - J \eta \rangle - \int_\Omega \langle \partial_t K \partial_t \psi, u \rangle\\
&=& -\int_\Omega \langle \partial_s \psi, (\partial_t J) u + \eta \rangle\\
&&\hspace{20pt}  - \int_\Omega \langle \partial_t \psi, -(\partial_s J) u - J \eta - \nu - \theta_\Sigma(\partial_t) \partial_z u \rangle .
\end{eqnarray*}
Here, since $(\partial_t J) u + \eta$ and $-(\partial_s J) u - J \eta - \nu - \theta_\Sigma(\partial_t) \partial_z u$ are in $L^r_{loc}$, the elliptic regularity of the usual Laplacian shows $u \in W^{1, r}_{loc}$ as desired. \qed

\vspace{2mm}
Iterating this proposition, we can obtain following heigher version regularity property.
We will prove in three steps to use the induction argument on integers $l$ and $k$.

\vspace{2mm}
\noindent
\textbf{Proposition 2.6}\hspace{2mm}
Fix a positive integer $l \ge 1$ and $k \in \{ 0, \ldots , l\}.$
Let $p > 3$ and assume $J \in W^{l, p}(\Omega; \R^{(2n + 1) \times (2n + 1)}), u \in L^p_{loc}(\Omega, \R^{2n + 1}), \eta, \nu \in W^{k, p}_{loc}(\Omega, \R^{2n + 1})$ satisfy the relation
\begin{eqnarray*}
&& \int_\Omega \langle \partial_s \varphi + J^T \partial_t \varphi, u \rangle = -\int_\Omega \langle \varphi, \eta + (\partial_t J)u \rangle,\\
&& \int_\Omega \langle \partial_t ((J^T)^2 + 1) \varphi, u \rangle = -\int_\Omega \langle \varphi, \nu + \theta_\Sigma(\partial_t) \partial_z u \rangle,\quad \mbox{and}\\
&& \partial_z u = \mbox{const}.
\end{eqnarray*}
for any $\varphi \in C^\infty_0(\Omega, \R^{2n + 1}).$
Then we have $u \in W^{k + 1, p}_{loc}(\Omega, \R^{2n + 1}).$

\vspace{1mm}
\noindent
{\bf Proof.}\hspace{2mm}
(Step 1).\hspace{1mm}
We first show that the claim holds when $k = 0$ and $l = 1$.

Note that for any $r \le p,$ we have $\eta, \nu \in W^{k, p}_{loc} \subset L^p_{loc} \subset L^r_{loc}.$ 
So, if $u \in L^q_{loc}$ and $\displaystyle \frac{p}{p - 1} < q < \frac{3p}{p - 3}$ then, by the Proposition 2.5 and the Sobolev's embedding theorem 2.4 , we have
\[ u \in W^{1, r}_{loc} \subset L^{q^\prime}_{loc},\]
where $\displaystyle q^\prime = \frac{3r}{3 - r},\ r = \frac{pq}{p + q}.$

Take $q_0$ so that $\displaystyle \frac{p}{p - 1} < q_0 \le p.$
Then $u \in L^p_{loc} \subset L^{q_0}_{loc}$ and by the previous argument $u \in L^{q_1}_{loc}$ where $\displaystyle q_1 = \frac{3r_0}{3 - r_0},\ r_0 = \frac{pq_0}{p + q_0}.$

If $\displaystyle q_1 < \frac{3p}{p - 3}$ then we have $u \in L^{q_2}_{loc},\ \displaystyle q_2 = \frac{3r_1}{3 - r_1},\ r_1 = \frac{pq_1}{p + q_1}.$
Therefore it follows by induction that $u \in W^{1, r_i}_{loc}$ for $j = 0,\ \ldots \ , m.$
With $j = m$ we have $\displaystyle q_m > \frac{3p}{p - 3}$ and hence $\displaystyle r_m = \frac{pq_m}{p + q_m} > 3,$ so $u$ is continuous.
Using Proposition 3.5 again with $q = \infty, r = p$ we obtain $u \in W^{1, p}_{loc}.$
This is the conclution case of the case $k = 0$ and $l = 1.$

\vspace{2mm}
\noindent
(Step 2).\hspace{1mm} Next we consider the case of $k = l = 1.$
Let $q, r > 1$ such that $\displaystyle \frac{1}{p} + \frac{1}{q} = \frac{1}{r}$.

Define $u^\prime = \partial_s u.$
Then $u \in W^{1, q}_{loc} \ \Rightarrow \ u^\prime \in L^q_{loc}.$
It is trivial that $\partial_z u^\prime = \partial_s \mbox{const} = 0.$
Moreover we can see there exist $\eta^\prime, \nu^\prime$ and $u^\prime$ also satisfies the assumption of Proposition 2.5

Indeed, 
\begin{eqnarray*}
\int \langle \partial_s \varphi + J^T \partial_t \varphi, u^\prime \rangle &=& -\int \langle \partial_s^2 \varphi + \partial_s J^T \partial_t \varphi, u \rangle\\
&=& -\int \langle \partial_s^2 \varphi + J^T \partial_t \partial_s \varphi, u \rangle - \int \langle (\partial_s J^T) \partial_t \varphi, u \rangle\\
&=& -\int \langle \varphi, \partial_s \eta + \partial_s (\partial_t J) u - \partial_t (\partial_s J)u \rangle\\
&=& -\int \langle \varphi, \eta^\prime + (\partial_t J)u^\prime \rangle,
\end{eqnarray*}
where $\eta^\prime = \partial_s \eta - (\partial_s J) \partial_t u$.
According to the H\"{o}lder inequality, $\partial_s \eta - (\partial_s J) \partial_t u$ is acutually in $L^r_{loc}.$
This is first condition for $u^\prime$ and $\eta^\prime.$

\vspace{2mm}
For the second condition,
\begin{eqnarray*}
\int \langle \partial_t K \varphi, u^\prime \rangle &=& -\int \langle \partial_t \partial_s K \varphi, u \rangle\\
&=& -\int \langle \partial_t K \partial_s \varphi, u \rangle - \int \langle \partial_t (\partial_s K) \varphi, u \rangle\\
&=& -\int \langle \varphi, \partial_s \nu + \partial_s \theta_\Sigma(\partial_t) \partial_z u + (\partial_s K)^T \partial_t u \rangle\\
&=& -\int \langle \varphi, \nu^\prime + \theta_\Sigma(\partial_t) \partial_z u^\prime \rangle,
\end{eqnarray*}
where $\nu^\prime = \partial_s \nu + (\partial_s \theta_\Sigma(\partial_t)) \partial_z u + (\partial_s K)^T \partial_t u \in L^r_{loc}.$

\vspace{2mm}
Hence, applying Proposition 2.5 to the mappings $u^\prime \in L^q_{loc}$ and $\eta^\prime, \nu \in L^r_{loc},$ we can see $u^\prime \in W^{1, r}_{loc}.$
To this end, due to Lemma 2.2, we have 
\[ \partial_t u = -J(\eta - \partial_s u) + \theta(\partial_t u) \xi = -J(\eta - u^\prime) + \theta_\Sigma(\partial_t) \partial_z u \in W^{1, r}_{loc}.\]
Therefore, $u \in W^{2, r}_{loc}.$
Using a sequence $\{ q_i\}$ and an argument which is same with Step 1, we conclude $u \in W^{2, p}_{loc}.$

\vspace{2mm}
\noindent
(Step 3).\hspace{1mm} We here see the result holds in general.
 
Assume, by induction, the proposition is true for some $k \ge 1$ and that $J \in W^{k + 1, p}$.
That is, the inductuin hypothesis is
\[ u \in L^p_{loc}\ \mbox{and}\ \eta,\ \nu \in W^{k, p}_{loc} \ \Rightarrow \ u \in W^{k + 1, p}_{loc}.\]
Let $\eta, \nu \in W^{k + 1, p}_{loc}.$
Define $u^\prime,\ \eta^\prime$ and $\nu^\prime$ as in Step 2 by
\begin{eqnarray*}
u^\prime &=& \partial_s u,\qquad \eta^\prime = \partial_s \eta - (\partial_s J) \partial_t u,\\
\nu^\prime &=& \partial_s \nu + (\partial_s \theta_\Sigma(\partial_t)) \partial_z u + (\partial_s K)^T \partial_t u.
\end{eqnarray*}
Apply the inducution hypothesis $u^\prime$ and $\eta^\prime, \nu^\prime  \in W^{k, p}_{loc},$ we have $u \in W^{k + 2, p}_{loc}.$
This prove the proposition. \qed

\vspace{2mm}
 Now, we are in a position to prove Theorem 1.
It suffices to prove the result in local coordinates on $\Sigma$ and in local coordinates on $M$ such that we have been used.
If $J$ is a $C^l$ almost complex strucutre on $P$, then the image of $J$ by a $C^l$ coordinate chart gives a $C^{l - 1}$ tensor field on $\R^{2n + 1}$ such that $J^2 = -1 + \zeta \partial/\partial z$, where the $1$-form $\zeta$ is a pullback of the contact form $\theta$ of $M$.

\vspace{2mm}
First, we assume that $k = 1$.
Then by the assumption of the theorem, we have $u \in W^{1, p}_{loc}(\Omega, \R^{2n + 1})$.
So we can use Proposition 2.6 with $k = 0$, if we take $\eta = (\theta_\Sigma(\partial_t) - \theta_\Sigma(\partial_s))\partial_z u ,\ \nu = 0,$ and $J$ replaced by $J \circ u \in W^{1, p}_{loc}(\Omega, \R^{2n + 1} \times \R^{2n + 1}).$

Hence $u \in W^{2, p}_{loc}(\Omega, \R^{2n + 1}),\ J \circ u \in W^{2, p}_{loc}(\Omega, \R^{2n + 1} \times \R^{2n + 1}).$
By induction on $k$ we can obtain that $u \in W^{k + 1, p}_{loc}(\Omega, \R^{2n + 1})$ for $k = 1,\ \ldots ,\ l - 1.$
Therefore $u \in  W^{l, p}_{loc}(\Omega, \R^{2n + 1}).$

\vspace{6mm}
\begin{center}
REFERENCES
\end{center}

\vspace{3mm}\noindent
[1] \, M. Audin, J. Lafontaine, \, {\it Holomorphic curves in symplectic geometry}, Progress in mathematics 117, 1994.

\vspace{3mm}\noindent
[2] \, E.Barletta, S. Dragomir, K. L. Duggal, \, {\it Foliations in Cauchy-Riemann Geometry}, Mathematical Surveys and Monographs, Vol140, American Mathematical Society, 2007. 

\vspace{3mm}\noindent
[3] \, D. E. Blair, \, {\it Riemaniann geometry of contact and symplectic manifolds}, Progress in mathematics 203, Birkh\"{a}user, 2002.

\vspace{3mm}\noindent
[4] \, S. Dragomir, G. Tomassini, \, {\it Differential Geometry and Analysis on CR Manifolds}, Progress in mathematics 246, Birkh\"{a}user, 2006.

\vspace{3mm}\noindent
[5] \, G. B. Folland, J. J. Kohn, \, {\it The Neumann Problem for the Cauchy-Riemann Complex}, Ann. of  Math. Studies, No 75, Princeton University Press, New Jersey, 1972.

\vspace{3mm}\noindent
[6] \, A. Floer, \, Morse theory for Lagrangian intersections, J. Diff. Geom. {\bf 28}(1988), 513-547.

\vspace{3mm}\noindent
[7] \, M. Gromov, \, Pseudo holomorphic curves in symplectic manifolds, Invent. Math. {\bf 82}(1985), 307-347.

\vspace{3mm}\noindent
[8] \, H\"{o}rmander, \, {\it Linear Partial Differential Operators}. Springer-Verlag, New York, 1963.

\vspace{3mm}\noindent
[9] \, D. McDuff, D. Salamon, \, {\it $J$-Holomorphic Curves and Quantum Cohomology}, University Lecture Series 6, A.M.S, 1995.

\vspace{3mm}\noindent
[10] McDuff D. , Salamon D. , {\it $J$-holomorphic curves and symplectic topology}. Colloquium publications vol. 52, American Mathematical Society, 2004. 

\vspace{3mm}\noindent
[11] \, R. Petit, \, Harmonic maps and strictly pseudoconvex CR manifolds, Comm. Ana. and Geom. {\bf 10}(2002) No 3, 575-610.

\vspace{3mm}\noindent
[12]\, T. Saotome,\, Removable singularity theorem of $J$-holomorphic mappings for strongly pseudo-convex manifolds,\, (preprint).

\vspace{3mm}\noindent
[13] \, R. Schoen, S.-T. Yau, \, {\it Lectures on Harmonic Maps}, Conference Proceedings and Lecture Notes in Geometry and Topology 2, International Press, Boston, 1997.

\vspace{3mm}\noindent
[14] \, N. Tanaka, \, {\it A Differential Geometric Study on Strongly Pseudo-Convex Manifolds}(Lectures in Mathematics, Department Mathematics, Kyoto University 9), \, Kinokuniya, Tokyo, Japan, 1975.

\end{document}